\documentclass[11pt]{amsart}
\usepackage{hyperref}
\usepackage[T1]{fontenc}

\usepackage{a4wide,amsmath,amssymb,latexsym,amsthm}
\usepackage{eucal}
\usepackage{comment}
\usepackage{graphicx}

\usepackage{color}

\newtheorem{theorem}{Theorem}[section]

\newtheorem{remark}[theorem]{Remark}
\newtheorem{lemma}[theorem]{Lemma}
\newtheorem{corollary}[theorem]{Corollary}

\graphicspath{{./Imgs/}}

\newcommand{\R}{\mathbb R}

\usepackage{pifont}


\begin{document}

\title[Coefficient inverse problem]{A unified Carleman approach to inverse problems for coupled parabolic and hyperbolic-ODE equations}


\author{Jone Apraiz}
\address{Department of  Mathematics, Faculty of Science and Technology, University of the Basque Country (EHU), 48940 Leioa, Spain.}
\email{jone.apraiz@ehu.eus}

\author{Rodrigo Lecaros}
\address{Universidad T\'ecnica Federico Santa Mar\'ia, Departamento de Matem\'atica, Casilla 110-V, Valpara\'iso, Chile.}
\email {rodrigo.lecaros@usm.cl}

\author{Sebasti\'an Zamorano}
\address{University of Deusto, Faculty of Engineering,  Avda Universidades 24, 48007 Bilbao (Bizkaia), Spain.}
\email{sebastian.zamorano@deusto.es}

\thanks{R. Lecaros was partially supported by FONDECYT(Chile) Grant NO: 1260574. J. Apraiz was supported by the project PID2024-158206NB-100, funded by MICIU$/$AEI$/$10.13039$/$501100011033 and FEDER, UE, the IMUS - María de Maeztu grant CEX2024-001517-M - Apoyo a Unidades de Excelencia María de Maeztu, funded by MICIU$/$AEI$/$10.13039$/$501100011033 and the Project of Consolidación Investigadora CNS2024-154725, also funded by MICIU$/$AEI$/$10.13039$/$501100011033. S. Zamorano was partially supported by the Basque Government through the consolidated research group ”CMRG-Computational Mathematics Research Group” under Grant IT1856-26 and by Grant PID2023-146872OB-I00-DyCMaMod of MICIU (Spain).}

\begin{abstract}
We study inverse coefficient problems for two evolution equations with relaxation: a Kelvin--Voigt viscoelastic equation and the Moore--Gibson--Thompson (MGT) equation arising in nonlinear acoustics. In both cases, we establish the Lipschitz stability for recovering a space-dependent coefficient from a single measurement.

Our approach exploits an exact factorization of each operator into a classical evolution operator (the heat operator for the viscoelastic system, or the wave operator for the MGT equation) coupled with a first-order differential operator in time. This makes it possible to combine global Carleman estimates with weighted ODE inequalities. For the Kelvin--Voigt model, the data consist of an interior dynamic observation and a Cauchy snapshot at a prescribed positive time. For the MGT equation, a single boundary trace is sufficient under a standard geometric condition. The resulting argument provides a common framework for coefficient inverse problems in evolution equations with relaxation effects.
\end{abstract}

\date{}
\maketitle

\section{Introduction}
\subsection{Motivation}
Inverse coefficient problems for evolution equations arise whenever a material parameter cannot be measured directly and must instead be inferred from the response of the system. Typical examples include damping, viscosity, and relaxation parameters in viscoelasticity, acoustics, and wave propagation. The parameter-to-observation map is nonlinear, and the resulting inverse problem is generally ill-posed, even when the corresponding forward problem is well posed. Establishing uniqueness and, more importantly, quantitative stability therefore requires analytical estimates that retain precise information on both the unknown coefficient and the available observation.

Global Carleman estimates, together with the Bukhgeim--Klibanov method introduced in \cite{BK}, constitute one of the principal tools for this purpose. The method reduces the nonlinear comparison of two coefficiens to a linear inverse source problem and then uses a weighted estimate to extract the coefficient from a suitable time trace. We refer to \cite{Klibanov2004,klibanov2013carleman,Yamamoto2025} for general accounts and to \cite{BellaYama} for hyperbolic systems. Related coefficient identification problems for parabolic equations can be found, for instance, in \cite{kravaris1985,kravaris1986,Lowe1994,Suzuki1980,Richter1981}.

The equations considered here are not themselves classical parabolic or hyperbolic equations. They contain damping and relaxation terms that produce a mixed PDE--ODE structure. The central observation of this paper is that each operator admits an exact decomposition into a heat or wave operator acting on a first-order differential expression in time, plus a lower-order term containing the unknown coefficient. This structure makes it possible to transfer classical Carleman estimates to the full relaxation model by combining them with a weighted ODE inequality. Related PDE--ODE viewpoints have proved useful for the Benjamin--Bona--Mahony equation, structurally damped equations, and models with memory; see \cite{yamamoto2003one,martin2013null,chaves2014null,Biccari2019,Biccari2020}.

For standard constant-coefficient realizations of both models, one branch of the characteristic roots approaches a finite relaxation rate as the spatial frequency tends to infinity. This loss of a uniform spectral gap is known to create difficulties in observability and controllability arguments for strongly damped and memory-type equations; see, among others, \cite{russell1986mathematical,leugering1989boundary,micu2001controllability, martin2013null,cerpa2018,lizama2019controllability}. Our proof does not rely on a spectral decomposition. Instead, it uses the PDE--ODE factorization at the level of the differential operator, allowing the parabolic and hyperbolic parts to be treated with weights adapted to their respective propagation mechanisms.

\subsection{The viscoelasticity inverse problem}\label{Visco_subsec}
The first model is the Kelvin--Voigt equation set in a bounded domain $\Omega \subset \mathbb{R}^N$, with $N=2$ or $N=3$, over a time interval $(0,T)$, with an additional space-dependent friction term:
\begin{align*}
\left\{
\begin{array}{ll}
u_{tt}-\Delta u-\Delta u_t+b(x)u_t=f
&\text{in }\Omega\times(0,T),\\
u=0&\text{on }\partial\Omega\times(0,T),\\
u(\cdot,0)=u_0,\quad u_t(\cdot,0)=u_1&\text{in }\Omega.
\end{array}
\right.
\end{align*}
Here $u$ is the displacement, $f$ is a source term, $u_0$ and $u_1$ stand for the initial displacement and velocity, respectively, $-\Delta u$ is the elastic restoring term, $-\Delta u_t$ represents Kelvin--Voigt damping, and $b(x)u_t$ models frictional damping. The simultaneous action of structural and frictional damping occurs in models for polymers, biological tissues, composites, and geomaterials; see \cite{Lagnese1989,chen1989proof,cavalcanti2004, liu1999semigroups}.

From a mathematical standpoint, the key observation is that the viscoelastic equation can be rewritten in the form
\begin{equation}\label{coupling_visco}
    H_bu:=H_0(Pu)+\widehat{b}(x)u_t=f,
\end{equation}
where
\begin{equation*}
H_0:=\partial_t-\Delta,\qquad P:=\partial_t+1,
\qquad \widehat b:=b-1.
\end{equation*}
Thus the strongly damped wave operator is expressed as a heat operator acting on $Pu$, coupled with the first-order relation $Pu=u_t+u$.

We keep the source and the initial data fixed and seek to determine $b$ from one experiment. Let $t_0\in(0,T)$ be a prescribed observation time and let $I=(t_0-\tau,t_0+\tau)\Subset(0,T)$, for $\tau >0$. The data are
\begin{equation*}
\left.(u_{tt}+u_t)\right|_{\omega\times I},
\qquad
\big(u(\cdot,t_0),u_t(\cdot,t_0)\big),
\end{equation*}
where $\omega\subset\Omega$ is a nonempty open set. Theorem \ref{mainresult2} gives a Lipschitz estimate for $b^1-b^2$ in $L^2(\Omega)$ in terms of the difference of these data. We assume that the two coefficients agree outside an interior subdomain $\Omega_0\Subset\Omega$ and that the reference velocity is bounded away from zero on $\Omega_0$ at time $t_0$. The localization of this non-degeneracy condition is compatible with the homogeneous Dirichlet boundary condition. The Cauchy snapshot is also structural: evaluating the difference equation at $t_0$ produces both $\Delta(u^1-u^2)(t_0)$ and $\Delta\partial_t(u^1-u^2)(t_0)$.

Inverse problems for viscoelasticity systems have also been studied in different contexts: identification of elastic and viscoelastic moduli in \cite{Lorenzi2006}, logarithmic stability results in \cite{Buhan2010}, and inverse source problems in \cite{Cheng2025}.

\subsection{The Moore-Gibson-Thompson inverse problem}\label{MGT_subsection}
The second model is the Moore--Gibson--Thompson equation posed in $\Omega \times (0,T)$. After normalizing the thermal relaxation parameter, it takes the form:
\begin{align*}
\left\{
\begin{array}{ll}
u_{ttt}+\alpha(x)u_{tt}-c^2\Delta u-b\Delta u_t=f
&\text{in }\Omega\times(0,T),\\
u=0&\text{on }\partial\Omega\times(0,T),\\
u(\cdot,0)=u_0,\quad u_t(\cdot,0)=u_1,\quad
u_{tt}(\cdot,0)=u_2&\text{in }\Omega,
\end{array}
\right.
\end{align*}
where $f$ is a source term, while $u_0, u_1$, and $u_2$ denote the initial displacement, velocity, and acceleration, respectively. This third-order equation models acoustic propagation in viscous and thermally relaxing media and is relevant, in particular, to high-intensity focused ultrasound. The constants $c,b>0$ describe the sound speed and diffusivity, respectively, whereas the space-dependent coefficient $\alpha$ accounts for damping and relaxation. We refer to \cite{Jordan2014,delloro2017on,kaltenbacher2011wellposedness} for the modeling and well-posedness background, and to \cite{kaltenbacher2012exponential,kaltenbacher2015mathematics, lasiecka2015moore,lasiecka2016moore} for further qualitative results.

A crucial structural observation is that the MGT operator can be decomposed as
\begin{equation*}\label{coupling_MGT}
    \mathcal L_\alpha u=L_0^{(b)}(P_a u)+\gamma_\alpha(x)u_{tt},
\end{equation*}
where
\begin{equation*}
L_0^{(b)}:=\partial_t^2-b\Delta,\qquad
P_a:=\partial_t+a,\qquad a:=\frac{c^2}{b},\qquad
\gamma_\alpha:=\alpha-a.
\end{equation*}
Thus the MGT operator is a wave operator with propagation speed $\sqrt{b}$ acting on $P_a u$, coupled with a
first-order ODE in time.

Here $b,c$, the source, and the initial data are fixed, and the objective is to recover $\alpha(x)$. On a relatively open part $\Gamma_0\subset\partial\Omega$, we consider the single boundary measurement
\begin{equation*}
\Lambda(\alpha)
:=\left.\partial_n(P_a u)\right|_{\Gamma_0\times(0,T)}.
\end{equation*}
Under the standard multiplier geometry and a sufficiently long observation time, Theorem \ref{mainresult} establishes a two-sided Lipschitz estimate between $\alpha^1-\alpha^2$ in $L^2(\Omega)$ and the time derivative of $\Lambda(\alpha^1)-\Lambda(\alpha^2)$. Equivalently, the inverse estimate uses the $H^1$-in-time norm of this one measured trace. The excitation condition is imposed on the prescribed initial acceleration $u_2$, which is exactly the factor multiplying the coefficient difference in the initial identity generated by the Bukhgeim--Klibanov argument. Previous work on MGT inverse problems includes \cite{arancibia2022boundary}; related controllability results can be found in \cite{lizama2019controllability,lizama2022boundary,
lizama2022controllability}.

\subsection{Main contributions}\label{Main_contri}
The principal contribution of this work is a unified Carleman strategy for two coefficient inverse problems whose principal dynamics are different but whose relaxation structure is analogous. The main results and the essential steps of their proofs are as follows.

\begin{itemize}
\item For the Kelvin--Voigt equation, we derive a Carleman estimate for the full operator (Theorem \ref{Carle_visco}) from a two-sided parabolic weight centered at $t_0$, the heat estimate of \cite{imanuvilov2009carleman}, and a weighted ODE inequality. A fixed-time trace estimate then links the Carleman-controlled variable to the coefficient. Differentiating the difference equation before making any change of unknown preserves the parabolic--ODE factorization and avoids the appearance of variable principal-order coefficients. This leads to the Lipschitz stability result of Theorem \ref{mainresult2}.

\item For the MGT equation, we combine the global wave Carleman estimate of \cite{baudouin2013global} with a weighted ODE estimate to obtain a Carleman inequality for the complete third-order operator (Theorem \ref{observability}). Applied to the time derivative of the difference of two solutions, this estimate recovers the coefficient from the initial acceleration and the single boundary trace $\Lambda(\alpha)$. Theorem
\ref{mainresult} yields both the inverse Lipschitz bound and the corresponding direct trace estimate.
\end{itemize}
\subsection{Outline of the article}
The remainder of the paper is organized as follows. Section \ref{Visco_section} introduces the Kelvin--Voigt inverse problem, develops the centered parabolic Carleman estimate and the fixed-time trace argument,
and proves the corresponding Lipschitz stability theorem. Section \ref{MGT_section} establishes the wave--ODE Carleman estimate and applies it to the MGT inverse problem. Finally, Section \ref{open} lists several open
questions suggested by the present analysis.

\subsection{Notation}
Throughout the paper, $\Omega\subset\R^N$, with $N=2$ or $N=3$, denotes a nonempty bounded domain with smooth boundary $\Gamma:=\partial\Omega$, and $T>0$ is a fixed final time. We denote by $n$ the outward unit normal vector to $\Gamma$ and by $d\sigma$ the surface measure on $\Gamma$. We write $\Omega_0\Subset\Omega$ to indicate that $\overline{\Omega_0}\subset \Omega$ is compactly contained.

The notation $\varphi_\lambda$ will always denote the Carleman weight associated with the problem under consideration: the parabolic weight \eqref{weight-calor} in Section \ref{Visco_section}, and the hyperbolic weight \eqref{1.8} in Section \ref{MGT_section}. Unless otherwise stated, $C$ denotes a generic positive constant that may change from line to line. Such a constant is independent of the functions being estimated and of the large Carleman parameter $s$; any relevant dependence on $\lambda$ or on the fixed a priori data will be specified when needed.

\section{The Kelvin--Voigt viscoelastic equation}\label{Visco_section}

Let $Q:=\Omega\times(0,T)$, $\Sigma:=\Gamma\times(0,T)$, $Q_\omega:=\omega\times(0,T)$, where $\omega\subset\Omega$ is a nonempty open set. We consider the Kelvin--Voigt equation
\begin{equation}\label{eq:1}
\left\{
\begin{array}{lll}
u_{tt}-\Delta u-\Delta u_t+b(x)u_t=f
&\text{in }Q,\\
u=0
&\text{on }\Sigma,\\
u(\cdot,0)=u_0,\qquad u_t(\cdot,0)=u_1
&\text{in }\Omega.
\end{array}
\right.
\end{equation}
Here $u$ denotes the displacement, $-\Delta u$ is the elastic restoring
term, $-\Delta u_t$ represents Kelvin--Voigt damping, and $b(x)u_t$
models spatially dependent frictional damping. The simultaneous presence of
the last two mechanisms is relevant in the description of polymers,
biological tissues, composite materials, and geomaterials; see, for instance,
\cite{Lagnese1989,chen1989proof,cavalcanti2004,liu1999semigroups}.

For $b\in L^\infty(\Omega)$ with $b\geq0$, the standard energy theory yields
well-posedness of \eqref{eq:1} for
\begin{equation*}
u_0\in H_0^1(\Omega),\qquad
u_1\in L^2(\Omega),\qquad
f\in L^2(0,T;H^{-1}(\Omega)).
\end{equation*}
In particular, the corresponding weak solution satisfies
\begin{equation*}
u\in C([0,T];H_0^1(\Omega))\cap C^1([0,T];L^2(\Omega)),
\qquad
u_t\in L^2(0,T;H_0^1(\Omega)).
\end{equation*}
The inverse argument below is first established for sufficiently regular
solutions and then extends, whenever appropriate, by the usual density
argument.

We keep the initial data and the source fixed and seek to determine the
coefficient $b$ from one experiment. The available data consist of an
interior time-dependent observation and a global Cauchy snapshot at a
prescribed time $t_0\in(0,T)$. The fixed-time information is essential:
the parabolic component of the Kelvin--Voigt equation does not allow a
uniform reconstruction of an arbitrary initial state from an $L^2$
observation.

Because the homogeneous Dirichlet condition forces the traces of sufficiently
regular velocities to vanish on $\Gamma$, a uniform non-degeneracy condition
cannot hold up to the boundary. We therefore use the standard assumption that
the unknown part of the coefficient is supported in an interior subdomain.

Fix $\Omega_0\Subset\Omega$, $M>0$, and a reference coefficient
$b_\ast\in W^{1,\infty}(\Omega)$ with $b_\ast\geq0$. We introduce the
admissible class
\begin{align}\label{asumavisco}
\mathcal B_M(b_\ast):=
\Big\{b\in W^{1,\infty}(\Omega):\; b\geq0,\ 
\|b\|_{W^{1,\infty}(\Omega)}\leq M,
b=b_\ast\ \text{a.e. in }\Omega\setminus\Omega_0\Big\}.
\end{align}

\subsection{A Carleman estimate centered at the observation time}

Fix $t_0\in(0,T)$ and choose
\begin{equation*}
0<\tau<\min\{t_0,T-t_0\},\qquad
t_-=t_0-\tau,\qquad t_+=t_0+\tau.
\end{equation*}
We write
\begin{equation*}
I=(t_-,t_+),\qquad Q_I=\Omega\times I,\qquad
\Sigma_I=\Gamma\times I,\qquad
Q_{\omega,I}=\omega\times I.
\end{equation*}
Choose $\eta\in C^2(\overline\Omega)$ satisfying
\begin{align*}
\eta>0\ \text{in }\Omega,\qquad
\eta=0\ \text{on }\Gamma,\qquad
|\nabla\eta|>0\ \text{in }\overline\Omega\setminus\omega_0,
\end{align*}
for some open set $\omega_0\Subset\omega$, where $\omega_0$ contains all critical points of $\eta$.
Define
\begin{align*}
\ell(t)&=(t-t_-)(t_+-t)=\tau^2-(t-t_0)^2,
\end{align*}
and, for $\lambda>0$, set
\begin{align}\label{weight-calor}
\alpha_\lambda(x,t)&=\frac{e^{\lambda\eta(x)}}{\ell(t)},\notag\\
\varphi_\lambda(x,t)
&=\frac{e^{\lambda\eta(x)}
-e^{2\lambda\|\eta\|_{L^\infty(\Omega)}}}{\ell(t)}.
\end{align}
The weight tends to $-\infty$ at both endpoints of $I$ and, for every
fixed $x$, attains its unique maximum at $t=t_0$. In particular, for
every fixed $\lambda>0$ there exists $c_\lambda>0$ such that
\begin{align}\label{strict_decay_visco}
\varphi_\lambda(x,t)
\leq\varphi_\lambda(x,t_0)-c_\lambda|t-t_0|^2,
\qquad (x,t)\in Q_I.
\end{align}

For a sufficiently regular function $Y$, we define
\begin{align*}
\mathcal I_{s,\lambda}(Y):=
\int_{Q_I}\Big((s\alpha_\lambda)^{-1}
  (|Y_t|^2+|\Delta Y|^2)
 +\lambda^2s\alpha_\lambda|\nabla Y|^2
 +\lambda^4(s\alpha_\lambda)^3|Y|^2\Big)
 e^{2s\varphi_\lambda}\,dxdt.
\end{align*}
The endpoint singularities are crucial: they eliminate the time-boundary
terms without imposing an initial or final condition. We use the standard
interior parabolic estimate of \cite{imanuvilov2009carleman}.

\begin{theorem}[\cite{imanuvilov2009carleman}]
\label{Carle_calor_visco}
There exists $\lambda_0>0$ such that, for every $\lambda\geq\lambda_0$,
there exist $s_0=s_0(\lambda)>0$ and $C>0$ for which
\begin{align}\label{eq:4}
\mathcal I_{s,\lambda}(Y)
\leq C\int_{Q_I}|H_0(Y)|^2e^{2s\varphi_\lambda}\,dxdt
 +C\int_{Q_{\omega,I}}\lambda^4(s\alpha_\lambda)^3
 |Y|^2e^{2s\varphi_\lambda}\,dxdt
\end{align}
for every $s\geq s_0$ and every
\begin{equation*}
Y\in H^1(I;L^2(\Omega))
 \cap L^2(I;H^2(\Omega)\cap H_0^1(\Omega)),
\qquad H_0:=\partial_t-\Delta.
\end{equation*}
No condition is imposed at $t=t_-$ or $t=t_+$. The constant $C$ is
independent of $s$ and $Y$.
\end{theorem}

We next record an ODE estimate adapted to a weight with an interior
maximum. Recall that $P=\partial_t+1$.

\begin{lemma}\label{centered_ODE_visco}
For every fixed $\lambda>0$, there exist $s_1>0$ and $C>0$ such that
every $v\in H^1(I;L^2(\Omega))$ satisfies
\begin{align}\label{centered_ODE_estimate_visco}
\int_{Q_I}|v|^2e^{2s\varphi_\lambda}\,dxdt
\leq\frac{C}{s}\int_{Q_I}|Pv|^2
e^{2s\varphi_\lambda}\,dxdt
+\frac{C}{\sqrt{s}}\int_\Omega|v(x,t_0)|^2
e^{2s\varphi_\lambda(x,t_0)}\,dx,
\quad s\geq s_1.
\end{align}
\end{lemma}

\begin{proof}
Set
$$
Y:=Pv=v_t+v.
$$
Since $v\in H^1(I;L^2(\Omega))$, the variation-of-constants formula holds in $L^2(\Omega)$:
\begin{equation}
\label{eq:variation-constants}
v(x,t)=e^{-(t-t_0)}v(x,t_0)
+
\int_{t_0}^{t}e^{-(t-r)}Y(x,r)dr.
\end{equation}
We estimate separately the two terms on the right-hand side.

Let
$$
D_\lambda(x)
:=
e^{2\lambda\|\eta\|_{L^\infty(\Omega)}}
-e^{\lambda\eta(x)}>0.
$$
Recalling that
$$
\ell(t)=\tau^2-(t-t_0)^2,
\qquad
\varphi_\lambda(x,t)
=
-\frac{D_\lambda(x)}{\ell(t)},
$$
we first derive a quantitative estimate for the difference of the weights. Let
$$
t_0\leq r\leq t<t_+,
\qquad
a:=t-t_0,
\qquad
b:=r-t_0.
$$

Then $0\leq b\leq a<\tau$, and
\begin{align*}
\varphi_\lambda(x,t)-\varphi_\lambda(x,r)
=-D_\lambda(x)
\left(\frac{1}{\tau^2-a^2}
-\frac{1}{\tau^2-b^2}
\right)=-D_\lambda(x)
\frac{a^2-b^2}
{(\tau^2-a^2)(\tau^2-b^2)}.
\end{align*}

Since
$$
(\tau^2-a^2)(\tau^2-b^2)\leq \tau^4
$$
and
$$
d_\lambda
:=
\min_{x\in\overline{\Omega}}D_\lambda(x)>0,
$$
we obtain
\begin{equation}
\label{eq:weight-difference-right}
\varphi_\lambda(x,t)-\varphi_\lambda(x,r)
\leq
-c_\lambda(a^2-b^2)
=
-c_\lambda(t-r)(t+r-2t_0),
\end{equation}
where $c_\lambda:=d_\lambda/\tau^4>0$.

The analogous estimate holds on the left half-interval. More precisely, if
$$
t_-<t\leq r\leq t_0,
\qquad
a:=t_0-t,
\qquad
b:=t_0-r,
$$
then $0\leq b\leq a<\tau$ and
\begin{equation}
\label{eq:weight-difference-left}
\varphi_\lambda(x,t)-\varphi_\lambda(x,r)
\leq
-c_\lambda(a^2-b^2).
\end{equation}

We now estimate the integral term in \eqref{eq:variation-constants}. Put
$$
F(x,r):=e^{s\varphi_\lambda(x,r)}Y(x,r).
$$
For $t\in(t_0,t_+)$, we have
\begin{align*}
e^{s\varphi_\lambda(x,t)}
\int_{t_0}^{t}e^{-(t-r)}Y(x,r)dr
=
\int_{t_0}^{t}
e^{-(t-r)}
e^{s(\varphi_\lambda(x,t)-\varphi_\lambda(x,r))}
F(x,r)\,dr.
\end{align*}
Thus, for each fixed $x$, the corresponding integral kernel is
$$
K_{s,x}^{+}(t,r)
:=
\mathbf{1}_{\{t_0\leq r\leq t\}}
e^{-(t-r)}
e^{s(\varphi_\lambda(x,t)-\varphi_\lambda(x,r))}.
$$

By \eqref{eq:weight-difference-right},
$$
\left|K_{s,x}^{+}(t,r)\right|
\leq
e^{-sc_\lambda(a^2-b^2)},
\qquad
a=t-t_0,\quad b=r-t_0.
$$

We next verify the two bounds required by Schur's test (see, e.g., \cite[Theorem 5.2]{MR517709}). Setting
$$
A=\sqrt{sc_\lambda}\,a,
\qquad
B=\sqrt{sc_\lambda}\,b,
$$
we obtain
\begin{align*}
\int_0^a e^{-sc_\lambda(a^2-b^2)}db
&=
\frac{1}{\sqrt{sc_\lambda}}
\int_0^A e^{-(A^2-B^2)}dB.
\end{align*}
The last integral is uniformly bounded with respect to $A\geq0$. Indeed, if $0\leq A\leq1$, then
$$
\int_0^A e^{-(A^2-B^2)}\,dB\leq A\leq1,
$$
whereas, if $A\geq1$, then $A^2-B^2=(A-B)(A+B)\geq A(A-B)$, and consequently
$$
\int_0^A e^{-(A^2-B^2)}dB
\leq
\int_0^A e^{-A(A-B)}dB
\leq \frac{1}{A}\leq1.
$$

It follows that
\begin{equation}
\label{eq:schur-first}
\sup_{0<a<\tau}
\int_0^a e^{-sc_\lambda(a^2-b^2)}db
\leq \frac{C}{\sqrt{s}}.
\end{equation}

For the second Schur bound, we similarly have
$$
\int_b^\tau e^{-sc_\lambda(a^2-b^2)}da
\leq
\frac{1}{\sqrt{sc_\lambda}}
\int_B^\infty e^{-(A^2-B^2)}dA.
$$
If $0\leq B\leq1$, then $\int_B^\infty e^{-(A^2-B^2)}dA \leq e\int_0^\infty e^{-A^2}dA$.
If $B\geq1$, then $A^2-B^2=(A-B)(A+B)\geq B(A-B)$, and hence
$$
\int_B^\infty e^{-(A^2-B^2)}dA
\leq
\int_B^\infty e^{-B(A-B)}dA
=
\frac{1}{B}\leq1.
$$

Therefore,
\begin{equation}
\label{eq:schur-second}
\sup_{0<b<\tau}
\int_b^\tau e^{-sc_\lambda(a^2-b^2)}da
\leq \frac{C}{\sqrt{s}}.
\end{equation}

It follows from \eqref{eq:schur-first}, \eqref{eq:schur-second}, and Schur's test that
$$
\left\|
e^{s\varphi_\lambda}
\int_{t_0}^{t}e^{-(t-r)}Y(r)dr
\right\|_{L^2(\Omega\times(t_0,t_+))}
\leq
\frac{C}{\sqrt{s}}
\left\|e^{s\varphi_\lambda}Y\right\|_{L^2(\Omega\times(t_0,t_+))}.
$$

For $t\in(t_-,t_0)$, we write
$$
\int_{t_0}^{t}e^{-(t-r)}Y(x,r)dr
=
-\int_t^{t_0}e^{r-t}Y(x,r)dr.
$$
Since $0\leq r-t\leq\tau$, the factor $e^{r-t}$ is bounded by $e^\tau$. Using \eqref{eq:weight-difference-left} and repeating the preceding Schur argument, we obtain
$$
\left\|
e^{s\varphi_\lambda}
\int_{t_0}^{t}e^{-(t-r)}Y(r)dr
\right\|_{L^2(\Omega\times(t_-,t_0))}
\leq
\frac{C}{\sqrt{s}}
\left\|e^{s\varphi_\lambda}Y\right\|_{L^2(\Omega\times(t_-,t_0))}.
$$

Combining the estimates on the two half-intervals and squaring gives
\begin{equation}
\label{eq:integral-term-estimate}
\int_{Q_I}
\left|
\int_{t_0}^{t}e^{-(t-r)}Y(x,r)dr
\right|^2
e^{2s\varphi_\lambda(x,t)}dxdt
\leq
\frac{C}{s}
\int_{Q_I}|Y|^2e^{2s\varphi_\lambda}dxdt.
\end{equation}

It remains to estimate the first term in \eqref{eq:variation-constants}. By the maximality of the weight at $t=t_0$, there exists, for fixed $\lambda>0$, a constant $c_\lambda>0$ such that
$$
\varphi_\lambda(x,t)
\leq
\varphi_\lambda(x,t_0)-c_\lambda|t-t_0|^2,
\qquad (x,t)\in Q_I.
$$
Moreover, $e^{-2(t-t_0)}\leq e^{2\tau}$ for $t\in I$.
Consequently,
\begin{align*}
\int_I
e^{-2(t-t_0)}
e^{2s\varphi_\lambda(x,t)}dt
&\leq
C e^{2s\varphi_\lambda(x,t_0)}
\int_I e^{-2sc_\lambda|t-t_0|^2}dt
\leq
C e^{2s\varphi_\lambda(x,t_0)}
\int_{\mathbb{R}}e^{-2sc_\lambda|t-t_0|^2}dt\\
&\leq
\frac{C}{\sqrt{s}}
e^{2s\varphi_\lambda(x,t_0)}.
\end{align*}
Hence
\begin{equation}
\label{eq:initial-term-estimate}
\int_{Q_I}
e^{-2(t-t_0)}
|v(x,t_0)|^2
e^{2s\varphi_\lambda(x,t)}dxdt
\leq
\frac{C}{\sqrt{s}}
\int_\Omega
|v(x,t_0)|^2e^{2s\varphi_\lambda(x,t_0)}dx.
\end{equation}

Finally, applying $|A+B|^2\leq 2|A|^2+2|B|^2$ to \eqref{eq:variation-constants}, and using \eqref{eq:integral-term-estimate} and \eqref{eq:initial-term-estimate}, we conclude that
$$
\int_{Q_I}|v|^2e^{2s\varphi_\lambda}dxdt
\leq
\frac{C}{s}
\int_{Q_I}|Pv|^2e^{2s\varphi_\lambda}dxdt
+
\frac{C}{\sqrt{s}}
\int_\Omega
|v(x,t_0)|^2e^{2s\varphi_\lambda(x,t_0)}dx
$$
for every $s\geq s_1$, where $s_1>0$ and $C>0$ may depend on $\lambda$, $\tau$, and the fixed weight, but are independent of $s$ and $v$.
\end{proof}

Next, we transfer \eqref{eq:4} to the viscoelastic operator
\begin{equation*}
H_bv:=v_{tt}-\Delta v-\Delta v_t+b(x)v_t.
\end{equation*}
Recall that $Pv=v_t+v$.

\begin{theorem}\label{Carle_visco}
Let $M>0$. There exists $\lambda_0>0$ such that, for every
$\lambda\geq\lambda_0$, there exist $s_0=s_0(\lambda,M)>0$ and $C>0$
such that
\begin{multline}\label{eq:5}
\mathcal I_{s,\lambda}(Pv)
\leq C\int_{Q_I}|H_bv|^2e^{2s\varphi_\lambda}\,dxdt
+C\int_{Q_{\omega,I}}\lambda^4(s\alpha_\lambda)^3
|Pv|^2e^{2s\varphi_\lambda}\,dxdt\\
+\frac{C}{\sqrt{s}}\int_\Omega|v(x,t_0)|^2
e^{2s\varphi_\lambda(x,t_0)}\,dx.
\end{multline}
for every $s\geq s_0$, every $b\in L^\infty(\Omega)$ satisfying
$\|b\|_{L^\infty(\Omega)}\leq M$, and every sufficiently regular $v$
such that $v=0$ on $\Sigma_I$.
The constant $C$ is uniform for $b$ in the indicated bounded set.
\end{theorem}

\begin{proof}
Set $Y:=Pv$ and $\widehat b:=b-1$. The factorization introduced in
\eqref{coupling_visco} gives
\begin{align}\label{factorization_visco}
H_bv=H_0(Y)+\widehat b\,v_t
    =H_0(Y)+\widehat b\,Y-\widehat b\,v.
\end{align}
Applying Theorem \ref{Carle_calor_visco} to $Y$ and using
\eqref{factorization_visco}, we obtain
\begin{multline}\label{intermediate_visco}
\mathcal I_{s,\lambda}(Y)
\leq C\int_{Q_I}|H_bv|^2e^{2s\varphi_\lambda}\,dxdt
 +C\int_{Q_I}\big(|Y|^2+|v|^2\big)
 e^{2s\varphi_\lambda}\,dxdt\\
 +C\int_{Q_{\omega,I}}\lambda^4(s\alpha_\lambda)^3
 |Y|^2e^{2s\varphi_\lambda}\,dxdt.
\end{multline}
Lemma \ref{centered_ODE_visco}, applied to $Pv=Y$, yields
\begin{align*}
\int_{Q_I}|v|^2e^{2s\varphi_\lambda}\,dxdt
\leq\frac{C}{s}\int_{Q_I}|Y|^2e^{2s\varphi_\lambda}\,dxdt
+\frac{C}{\sqrt{s}}\int_\Omega|v(x,t_0)|^2
e^{2s\varphi_\lambda(x,t_0)}\,dx.
\end{align*}
Furthermore, $\alpha_\lambda$ has a positive lower bound in $Q_I$.
Applying Lemma~\ref{centered_ODE_visco} to $v$, the term $\int_{Q_I} |v|^2 e^{2s\varphi_\lambda}\,dxdt$ is bounded by
\begin{equation*}
\frac{C}{s}\int_{Q_I}|Y|^2e^{2s\varphi_\lambda}\,dxdt
+\frac{C}{\sqrt{s}}\int_\Omega|v(x,t_0)|^2 e^{2s\varphi_\lambda(x,t_0)}\,dx.
\end{equation*}
Since $s\alpha_\lambda \ge c > 0$, the integral $\int_{Q_I}\big(|Y|^2+|v|^2\big) e^{2s\varphi_\lambda}\,dxdt$ in \eqref{intermediate_visco} is dominated by $\frac{C}{s} \mathcal I_{s,\lambda}(Y)$ plus the boundary term at $t_0$. Thus, choosing $s_0$ sufficiently large, the term $\frac{C}{s}\mathcal I_{s,\lambda}(Y)$ is absorbed by the left-hand side of \eqref{intermediate_visco}, leaving the snapshot trace at $t=t_0$ explicitly on the right-hand side. This yields \eqref{eq:5}.
\end{proof}

The recovery argument also uses the following weighted trace estimate.

\begin{lemma}\label{interior_trace_visco}
For every fixed $\lambda\geq\lambda_0$, there exist $s_1>0$ and $C>0$
such that
\begin{align}\label{trace_estimate_visco}
\int_\Omega|Y(x,t_0)|^2e^{2s\varphi_\lambda(x,t_0)}\,dx
\leq C\mathcal I_{s,\lambda}(Y)
\end{align}
for all $s\geq s_1$ and every $Y$ in the space of
Theorem \ref{Carle_calor_visco}.
\end{lemma}

\begin{proof}
Fix $\lambda\geq\lambda_0$ and first assume that $Y$ is smooth. We
divide the argument into three steps.

\smallskip
\noindent
\emph{Step 1: integration by parts in time.}
Set
\begin{equation*}
F(t):=\int_\Omega |Y(x,t)|^2
e^{2s\varphi_\lambda(x,t)}\,dx.
\end{equation*}
The endpoint singularities in \eqref{weight-calor} imply that
$F(t)\to0$ as $t\to t_-^+$ and as $t\to t_+^-$. Moreover,
\begin{equation*}
F'(t)=2\int_\Omega YY_t e^{2s\varphi_\lambda}\,dx
+2s\int_\Omega(\partial_t\varphi_\lambda)|Y|^2
e^{2s\varphi_\lambda}\,dx.
\end{equation*}
Integrating this identity on $(t_-,t_0)$ gives
\begin{align}\label{trace_identity_detailed}
\int_\Omega|Y(x,t_0)|^2e^{2s\varphi_\lambda(x,t_0)}\,dx
=2\int_{t_-}^{t_0}\!\!\int_\Omega YY_t
e^{2s\varphi_\lambda}\,dxdt
+2s\int_{t_-}^{t_0}\!\!\int_\Omega
(\partial_t\varphi_\lambda)|Y|^2
e^{2s\varphi_\lambda}\,dxdt.
\end{align}

\smallskip
\noindent
\emph{Step 2: estimates for the two terms in
\eqref{trace_identity_detailed}.}
Direct differentiation of \eqref{weight-calor}, together with the fact
that $|\ell'|$ is bounded, gives
\begin{align}\label{weight_derivative_detailed}
|\partial_t\varphi_\lambda|
\leq C\alpha_\lambda^2
\qquad\text{in }Q_I.
\end{align}
Using
\begin{equation*}
2|YY_t|
\leq(s\alpha_\lambda)^{-1}|Y_t|^2
+s\alpha_\lambda|Y|^2
\end{equation*}
and \eqref{weight_derivative_detailed} in
\eqref{trace_identity_detailed}, we deduce
\begin{align}\label{trace_intermediate_detailed}
\int_\Omega|Y(x,t_0)|^2e^{2s\varphi_\lambda(x,t_0)}\,dx
\leq C\int_{Q_I}\left((s\alpha_\lambda)^{-1}|Y_t|^2
+(s\alpha_\lambda+s\alpha_\lambda^2)|Y|^2\right)
e^{2s\varphi_\lambda}\,dxdt.
\end{align}

It remains to explain why the last two terms are controlled by the
zero-order part of $\mathcal I_{s,\lambda}$. There exists
$\alpha_\ast=\alpha_\ast(\lambda)>0$ such that
\begin{equation*}
\alpha_\lambda(x,t)\geq\alpha_\ast
\qquad\text{in }Q_I.
\end{equation*}
Consequently, for $s\geq1$,
\begin{align*}
s\alpha_\lambda
\leq\frac{1}{s^2\alpha_\ast^2}
(s\alpha_\lambda)^3,
\quad
s\alpha_\lambda^2
\leq\frac{1}{s^2\alpha_\ast}
(s\alpha_\lambda)^3.
\end{align*}
As $\lambda$ is fixed, these inequalities and
\eqref{trace_intermediate_detailed} give
\begin{equation*}
\int_\Omega|Y(x,t_0)|^2e^{2s\varphi_\lambda(x,t_0)}\,dx
\leq C\mathcal I_{s,\lambda}(Y)
\end{equation*}
for all sufficiently large $s$.

\smallskip
\noindent
\emph{Step 3: passage to the stated regularity class.}
Let
\begin{equation*}
\mathcal X:=H^1(I;L^2(\Omega))
\cap L^2(I;H^2(\Omega)\cap H_0^1(\Omega)).
\end{equation*}
For $Y\in\mathcal X$, choose smooth functions $Y_m$, satisfying the
same homogeneous boundary condition, such that $Y_m\to Y$ in
$\mathcal X$. The trace map
\begin{equation*}
H^1(I;L^2(\Omega))\longrightarrow C(\overline I;L^2(\Omega))
\end{equation*}
is continuous; hence $Y_m(\cdot,t_0)\to Y(\cdot,t_0)$ in
$L^2(\Omega)$.
Furthermore, for fixed $s$ and $\lambda$, all coefficients
\begin{equation*}
e^{2s\varphi_\lambda}(s\alpha_\lambda)^{-1},\qquad
e^{2s\varphi_\lambda}s\alpha_\lambda,
\qquad
e^{2s\varphi_\lambda}(s\alpha_\lambda)^3
\end{equation*}
are bounded in $Q_I$. Near $t=t_-$ and $t=t_+$, this follows from the
exponential decay displayed in Step~1, which dominates every power of
$\ell(t)^{-1}$. Therefore
$\mathcal I_{s,\lambda}(Y_m-Y)\to0$, and we may pass to the limit in
the estimate for $Y_m$. This proves \eqref{trace_estimate_visco} for
every $Y\in\mathcal X$.
\end{proof}

\subsection{Lipschitz stability for the damping coefficient}

The fixed-time data are the Cauchy snapshot
$(u(\cdot,t_0),u_t(\cdot,t_0))$. Both components have a structural
role: the equation at $t=t_0$ contains both $\Delta z(\cdot,t_0)$ and
$\Delta z_t(\cdot,t_0)$, where $z$ is the difference of two solutions.
The dynamic observation is $Pz_t=z_{tt}+z_t$, the parabolic variable
controlled by Theorem \ref{Carle_visco}.

\begin{theorem}\label{mainresult2}
Let $\Omega_0\Subset\Omega$ and let $\mathcal B_M(b_\ast)$ be defined by
\eqref{asumavisco}. For $k=1,2$, let $b^k\in\mathcal B_M(b_\ast)$ and
let $u^k$ be the solution of \eqref{eq:1} corresponding to $b=b^k$, for $k=1,2$, with the same source and initial data. Assume that, for some $K>0$,
\begin{align}\label{regsolvisco}
u^k\in H^3(0,T;H^2(\Omega)\cap H_0^1(\Omega)),
\qquad
\|u^k\|_{H^3(0,T;H^2(\Omega))}\leq K,
\quad k=1,2.
\end{align}
Suppose, in addition, that the reference solution satisfies
\begin{align}\label{u_t^2}
|u_t^2(x,t_0)|\geq c_0
\qquad\text{for a.e. }x\in\Omega_0
\end{align}
for some $c_0>0$. Then there exists a constant
\begin{equation*}
C=C(\Omega,T,t_0,\tau,\omega,\Omega_0,M,K,c_0)>0
\end{equation*}
such that
\begin{multline}\label{stabivisco}
\|b^1-b^2\|_{L^2(\Omega)}
\leq C\Big(
\|\partial_t^2(u^1-u^2)+\partial_t(u^1-u^2)\|_{L^2(Q_{\omega,I})}\\
+\|(u^1-u^2)(\cdot,t_0)\|_{H^2(\Omega)}
+\|\partial_t(u^1-u^2)(\cdot,t_0)\|_{H^2(\Omega)}
\Big).
\end{multline}
Consequently, the first term on the right-hand side may be replaced by
$\|u^1-u^2\|_{H^2(0,T;L^2(\omega))}$.
\end{theorem}

\begin{remark}
The condition $b^1=b^2=b_\ast$ in $\Omega\setminus\Omega_0$ is the
usual known-near-the-boundary assumption in coefficient inverse problems.
It is compatible with \eqref{u_t^2}, whereas a uniform lower bound for
$u_t^2(\cdot,t_0)$ on the whole of $\Omega$ would conflict with the
homogeneous Dirichlet boundary condition for regular solutions.
\end{remark}

\begin{remark}
Assumption \eqref{regsolvisco} is used to justify the time differentiation
of the difference equation and to obtain a uniform bound for
$u_{tt}^2$ in $Q$. Since $N\leq3$, this bound follows from
$H^2(\Omega)\hookrightarrow L^\infty(\Omega)$ and the one-dimensional
embedding $H^1(0,T;H^2(\Omega))\hookrightarrow
C([0,T];H^2(\Omega))$.
\end{remark}

\begin{remark}
The use of the Cauchy pair at $t_0$ is not redundant. With only
$z(\cdot,t_0)$, the term $\Delta z_t(\cdot,t_0)$ in
\eqref{coefficient_identity_visco} below would remain uncontrolled.
Removing the velocity snapshot would require additional structure or a
substantially stronger observation theorem.
\end{remark}

\begin{proof}[Proof of Theorem \ref{mainresult2}]
Set
\begin{equation*}
z:=u^1-u^2,\qquad q:=b^2-b^1,\qquad R:=u_t^2.
\end{equation*}
Subtracting the two equations gives
\begin{equation}\label{difference_visco}
\left\{
\begin{array}{lll}
H_{b^1}z=q(x)R(x,t)&\text{in }Q,\\
z=0&\text{on }\Sigma,\\
z(\cdot,0)=z_t(\cdot,0)=0&\text{in }\Omega.
\end{array}
\right.
\end{equation}
Since $q$ is independent of time, the function $v:=z_t$ satisfies
\begin{equation}\label{differentiated_visco}
\left\{
\begin{array}{lll}
H_{b^1}v=q(x)u_{tt}^2(x,t)&\text{in }Q_I,\\
v=0&\text{on }\Sigma_I.
\end{array}
\right.
\end{equation}
On $Q_I$, this equation is used without an endpoint condition. Set
$Y=Pv=v_t+v=z_{tt}+z_t$. Evaluating \eqref{difference_visco} at
$t=t_0$ gives the exact coefficient identity
\begin{align}\label{coefficient_identity_visco}
qR(\cdot,t_0)
=Y(\cdot,t_0)-\Delta z(\cdot,t_0)
-\Delta v(\cdot,t_0)+(b^1-1)v(\cdot,t_0).
\end{align}

Fix $\lambda\geq\lambda_0$ and take $s$ sufficiently large. Applying
Theorem \ref{Carle_visco} to \eqref{differentiated_visco}, and then
Lemma \ref{interior_trace_visco} with $Y=Pv$, yields
\begin{multline}\label{weighted_recovery_visco}
\int_\Omega|Y(x,t_0)|^2e^{2s\varphi_\lambda(x,t_0)}\,dx
\leq C\int_{Q_I}|q(x)u_{tt}^2(x,t)|^2
e^{2s\varphi_\lambda}\,dxdt\\
+C\int_{Q_{\omega,I}}\lambda^4(s\alpha_\lambda)^3
|Y|^2e^{2s\varphi_\lambda}\,dxdt
+\frac{C}{\sqrt{s}}\int_\Omega|v(x,t_0)|^2
e^{2s\varphi_\lambda(x,t_0)}\,dx.
\end{multline}
By \eqref{asumavisco}, $q=0$ a.e. in
$\Omega\setminus\Omega_0$.
The regularity assumption \eqref{regsolvisco} and the embeddings recalled
above give $\|u_{tt}^2\|_{L^\infty(Q)}\leq C K$. In order to estimate the first term in the right-hand side of \eqref{weighted_recovery_visco}, since $I=(t_0-\tau,t_0+\tau)$, set $r=\sqrt{2sc_{\lambda}}(t-t_0)$, we obtain
\begin{align*}
\int_{I} e^{-2sc_\lambda|t-t_0|^2}dt\leq \frac{1}{\sqrt{2sc_\lambda}}\int_{\R}e^{-r^2}dr=\sqrt{\frac{\pi}{2c_\lambda}}\frac{1}{\sqrt{s}}.
\end{align*}
Therefore, 
\begin{align}\label{source_absorption_visco}
\int_{Q_I}|q(x)u_{tt}^2(x,t)|^2
e^{2s\varphi_\lambda(x,t)}\,dxdt
\leq\frac{CK^2}{\sqrt{s}}
\int_{\Omega_0}|q(x)|^2
e^{2s\varphi_\lambda(x,t_0)}\,dx,
\end{align}
where $C$ is a positive constant independent of $s$.

On the other hand, \eqref{coefficient_identity_visco},
$\|b^1\|_{L^\infty}\leq M$, and \eqref{u_t^2} imply
\begin{multline}\label{nondegenerate_lower_visco}
c_0^2\int_{\Omega_0}|q|^2
e^{2s\varphi_\lambda(\cdot,t_0)}\,dx
\leq C\int_\Omega|Y(\cdot,t_0)|^2
e^{2s\varphi_\lambda(\cdot,t_0)}\,dx\\
+C\int_\Omega\big(|\Delta z(\cdot,t_0)|^2
+|\Delta v(\cdot,t_0)|^2+|v(\cdot,t_0)|^2\big)
e^{2s\varphi_\lambda(\cdot,t_0)}\,dx.
\end{multline}
Choosing $s$ so large that the right-hand side of
\eqref{source_absorption_visco} can be absorbed into
\eqref{nondegenerate_lower_visco}, we obtain
\begin{multline}\label{final_weighted_visco}
\int_{\Omega_0}|q|^2e^{2s\varphi_\lambda(\cdot,t_0)}\,dx
\leq C\int_{Q_{\omega,I}}\lambda^4(s\alpha_\lambda)^3
|z_{tt}+z_t|^2e^{2s\varphi_\lambda}\,dxdt\\
+C\big(\|z(\cdot,t_0)\|_{H^2(\Omega)}^2
+\|z_t(\cdot,t_0)\|_{H^2(\Omega)}^2\big).
\end{multline}
Finally, $\lambda$ and $s$ are now fixed. The weight on the left-hand side
of \eqref{final_weighted_visco} is bounded from below by a positive
constant, while the weighted factor on the right-hand side is bounded
from above on $Q_I$. Since $q$ vanishes outside $\Omega_0$, taking
square roots gives \eqref{stabivisco}.
\end{proof}

\begin{corollary}\label{uniqueness_visco}
Under the assumptions of Theorem \ref{mainresult2}, if
\begin{equation*}
\partial_t^2u^1+\partial_tu^1
=\partial_t^2u^2+\partial_tu^2
\qquad\text{in }Q_{\omega,I},
\end{equation*}
and
\begin{equation*}
(u^1,u_t^1)(\cdot,t_0)=(u^2,u_t^2)(\cdot,t_0)
\qquad\text{in }\Omega,
\end{equation*}
then $b^1=b^2$ a.e. in $\Omega$.
\end{corollary}


\section{The Moore--Gibson--Thompson equation}\label{MGT_section}

The Moore--Gibson--Thompson (MGT) equation is a third-order model for
acoustic propagation in viscous and thermally relaxing media. It arises when
Fourier's law is replaced by the Cattaneo--Vernotte law, thereby introducing a
finite thermal relaxation time. The model is relevant, for example, in
high-intensity focused ultrasound and nonlinear acoustics; see
\cite{Jordan2014,delloro2017on,kaltenbacher2011wellposedness}. Further results
on well-posedness, long-time behavior, and controllability can be found in
\cite{kaltenbacher2012exponential,kaltenbacher2015mathematics,
lasiecka2015moore,lasiecka2016moore,lizama2019controllability,
lizama2022boundary,lizama2022controllability}.

We consider
\begin{equation}\label{eqmain2}
\left\{
\begin{array}{lll}
u_{ttt}+\alpha(x)u_{tt}-c^2\Delta u-b\Delta u_t=f
&\text{in }Q,\\
u=0
&\text{on }\Sigma,\\
u(\cdot,0)=u_0,\quad u_t(\cdot,0)=u_1,\quad
u_{tt}(\cdot,0)=u_2
&\text{in }\Omega,
\end{array}
\right.
\end{equation}
where $b,c>0$ are known constants and $\alpha=\alpha(x)$ is the
coefficient to be recovered. Recalling 
\begin{equation*}
a:=\frac{c^2}{b},\qquad
\gamma_\alpha(x):=\alpha(x)-a,
\end{equation*}
define
\begin{equation*}
L_0^{(b)}:=\partial_t^2-b\Delta,
\qquad
P_a:=\partial_t+a.
\end{equation*}
A direct computation gives the exact factorization
\begin{align}\label{MGT_factorization}
\mathcal L_\alpha u
&:=u_{ttt}+\alpha u_{tt}-c^2\Delta u-b\Delta u_t\notag\\
&=L_0^{(b)}(P_a u)+\gamma_\alpha(x)u_{tt}.
\end{align}
This is the hyperbolic--ODE structure used throughout the section.

For instance, if $\alpha\in L^\infty(\Omega)$ and
\begin{align*}\label{dat}
u_0\in H^2(\Omega)\cap H_0^1(\Omega),\qquad
u_1\in H_0^1(\Omega),\qquad
u_2\in L^2(\Omega),\qquad
f\in L^2(Q),
\end{align*}
then the standard MGT well-posedness theory gives
\begin{equation*}
(u,u_t,u_{tt})\in
C\big([0,T];(H^2(\Omega)\cap H_0^1(\Omega))
\times H_0^1(\Omega)\times L^2(\Omega)\big);
\end{equation*}
see \cite{kaltenbacher2012exponential}. Additional regularity will be
specified in the inverse result.

Let $\Gamma_0\subset\Gamma$ be relatively open and set
$\Sigma_0:=\Gamma_0\times(0,T)$. The observation operator is
\begin{align}\label{MGT_measurement}
\Lambda(\alpha):=
\left.\partial_n(P_a u)\right|_{\Sigma_0}
=\left.\partial_n\left(u_t+\frac{c^2}{b}u\right)\right|_{\Sigma_0}.
\end{align}
The inverse problem is to determine $\alpha$ from this single time-dependent
boundary measurement, with $b,c,f$ and the initial data fixed.

Fix $M>a$ and introduce the admissible set
\begin{align*}
\mathcal B_M:=
\left\{\alpha\in L^\infty(\Omega):
a\leq\alpha(x)\leq M\ \text{for a.e. }x\in\Omega\right\}.
\end{align*}
Thus $0\leq\gamma_\alpha\leq M-a$ uniformly in this class.

\subsection{The wave and MGT Carleman estimates}

We impose the standard multiplier geometry: there exists
$x_0\in\mathbb R^N\setminus\overline\Omega$ such that
\begin{align}\label{asumb}
\Gamma_0\supset
\Gamma_+(x_0):=
\{x\in\Gamma:(x-x_0)\cdot n(x)\geq0\}.
\end{align}
Set
\begin{equation*}
R_0:=\sup_{x\in\Omega}|x-x_0|.
\end{equation*}
Because the propagation speed of $L_0^{(b)}$ is $\sqrt b$, the corresponding
time condition is
\begin{align}\label{asumc}
T>\frac{R_0}{\sqrt b}.
\end{align}
Choose $\beta\in(0,1)$ such that
\begin{align}\label{TC22}
R_0<\beta\sqrt b\,T.
\end{align}
For a sufficiently large constant $C_0>0$, define
\begin{align}\label{1.9}
\psi_b(x,t):=|x-x_0|^2-\beta b t^2+C_0,
&\quad \psi_b\geq1\quad\text{in }\overline Q,\\
\label{1.8}
\varphi_\lambda(x,t):=e^{\lambda\psi_b(x,t)},
&\quad\lambda>0.
\end{align}
The scaling in \eqref{1.9} follows by applying the unit-speed wave estimate
after the change of variable $\theta=\sqrt b\,t$.

For fixed $\lambda>0$ and $s>0$, let
\begin{multline*}
\mathcal J_s(Z):=
\sqrt{s}\int_\Omega e^{2s\varphi_\lambda(x,0)}
 |Z_t(x,0)|^2\,dx
+s\int_Qe^{2s\varphi_\lambda}\\
 \big(|Z_t|^2+b|\nabla Z|^2\big)\,dxdt
+s^3\int_Qe^{2s\varphi_\lambda}|Z|^2\,dxdt.
\end{multline*}

The following estimate is the constant-speed version of
\cite[Theorem 1.2]{baudouin2013global}, obtained by the preceding time
rescaling.

\begin{theorem}\label{CW2}
Assume \eqref{asumb}--\eqref{TC22}. There exists $\lambda_0>0$ such
that, for every fixed $\lambda\geq\lambda_0$, there exist $s_0>0$ and
$C>0$ satisfying
\begin{align}\label{1.13}
\mathcal J_s(Z)
\leq C\int_Qe^{2s\varphi_\lambda}
 |L_0^{(b)}Z|^2\,dxdt
+Cs\int_{\Sigma_0}e^{2s\varphi_\lambda}
 |\partial_nZ|^2\,d\sigma dt
\end{align}
for all $s\geq s_0$ and all admissible $Z$ such that
\begin{equation*}
Z=0\ \text{on }\Sigma,\qquad Z(\cdot,0)=0.
\end{equation*}

In particular, the estimate holds for
$Z\in L^2(0,T;H_0^1(\Omega))$ with $Z(\cdot,0)=0$,
whenever $L_0^{(b)}Z\in L^2(Q)$ and
$\partial_nZ\in L^2(\Sigma_0)$.
\end{theorem}

As in the previous section, an ODE-type Carleman inequality is needed in order to transfer \eqref{1.13} to the MGT operator.  It is a space-dependent version of the Volterra estimate in \cite[Lemma 3.1]{klibanov2013carleman}, adapted here to the first-order operators arising in our model.

\begin{lemma}\label{lemma_ODE}
Let $a\geq 0$ and let $\psi\in C^1(\overline{\Omega}\times[0,T))$ satisfy
\begin{align}\label{decreasing_weight}
\partial_t\psi(x,t)\leq -2\kappa t,
\qquad (x,t)\in\overline{\Omega}\times[0,T),
\end{align}
for some constant $\kappa>0$. Given $z\in L^2(\Omega\times(0,T))$, let $y\in H^1(0,T;L^2(\Omega))$ be the solution of
\begin{align}\label{first_order_problem}
\left\{
\begin{array}{ll}
P_a y=z, & \Omega\times(0,T),\\
y(\cdot,0)=0, & \Omega.
\end{array}
\right.
\end{align}
Then, for every $s>0$,
\begin{align}\label{desigualdadLemaEDO}
\int_0^T\int_\Omega |y|^2e^{2s\psi}\,dxdt
\leq \frac{1}{4s\kappa}
\int_0^T\int_\Omega |z|^2e^{2s\psi}\,dxdt.
\end{align}
\end{lemma}

\begin{proof}
For almost every $x\in\Omega$, the solution of \eqref{first_order_problem} is given by
\begin{align*}
y(x,t)=\int_0^t e^{-a(t-\tau)}z(x,\tau)\,d\tau.
\end{align*}
Since $a\geq0$, the Cauchy--Schwarz inequality yields
\begin{align}\label{ODE_pointwise_bound}
|y(x,t)|^2
&\leq t\int_0^t e^{-2a(t-\tau)}|z(x,\tau)|^2\,d\tau\notag\\
&\leq t\int_0^t |z(x,\tau)|^2\,d\tau.
\end{align}
Moreover, by \eqref{decreasing_weight},
\begin{align*}
t e^{2s\psi(x,t)}
\leq -\frac{1}{4s\kappa}\partial_t\left(e^{2s\psi(x,t)}\right).
\end{align*}
Consequently, for every $\tau\in(0,T)$ and $T'\in(\tau,T)$,
\begin{align}\label{weighted_time_integral}
\int_\tau^{T'} t e^{2s\psi(x,t)}\,dt
&\leq \frac{1}{4s\kappa}
\left(e^{2s\psi(x,\tau)}-e^{2s\psi(x,T')}\right)\notag\\
&\leq \frac{1}{4s\kappa}e^{2s\psi(x,\tau)}.
\end{align}
Letting $T'\to T$ in \eqref{weighted_time_integral}, we obtain the same estimate with $T$ as the upper integration limit. This argument also covers weights that are singular at $t=T$.
Multiplying \eqref{ODE_pointwise_bound} by $e^{2s\psi(x,t)}$, integrating over $\Omega\times(0,T)$, and applying Tonelli's theorem and \eqref{weighted_time_integral}, we obtain
\begin{align*}
\int_0^T\int_\Omega |y|^2e^{2s\psi}\,dxdt
&\leq \int_\Omega\int_0^T |z(x,\tau)|^2
\left(\int_\tau^T t e^{2s\psi(x,t)}\,dt\right)d\tau dx\\
&\leq \frac{1}{4s\kappa}
\int_0^T\int_\Omega |z|^2e^{2s\psi}\,dxdt,
\end{align*}
which proves \eqref{desigualdadLemaEDO}.
\end{proof}

\begin{remark}\label{remark_ODE_weights}
The hyperbolic weight satisfies the hypothesis of Lemma \ref{lemma_ODE}. Indeed,
\begin{align}\label{MGT_weight_decay}
\partial_t\varphi_\lambda(x,t)
=-2\beta b\lambda t\,\varphi_\lambda(x,t)
\leq-2\beta b\lambda e^\lambda t.
\end{align}
Thus \eqref{decreasing_weight} holds with
$\kappa=\beta b\lambda e^\lambda$.
\end{remark}

We now transfer the wave estimate to the full MGT operator.

\begin{theorem}\label{observability}
Assume \eqref{asumb}--\eqref{TC22} and let $M>a$. There exists
$\lambda_0>0$ such that, for every fixed $\lambda\geq\lambda_0$, there
exist $s_0=s_0(M)>0$ and $C>0$ for which
\begin{align}\label{1.133}
\mathcal J_s(P_a y)
\leq C\int_Qe^{2s\varphi_\lambda}
 |\mathcal L_\alpha y|^2\,dxdt
+Cs\int_{\Sigma_0}e^{2s\varphi_\lambda}
 |\partial_n(P_a y)|^2\,d\sigma dt
\end{align}
for all $s\geq s_0$, all $\alpha\in\mathcal B_M$, and every
sufficiently regular $y$ satisfying
\begin{equation*}
y=0\ \text{on }\Sigma,\qquad
y(\cdot,0)=y_t(\cdot,0)=0.
\end{equation*}
The constant $C$ is uniform with respect to $\alpha\in\mathcal B_M$.
\end{theorem}

\begin{proof}
Set
\begin{equation*}
Z:=P_a y,\qquad \gamma:=\gamma_\alpha.
\end{equation*}
Since $y(\cdot,0)=y_t(\cdot,0)=0$, we have $Z(\cdot,0)=0$. Moreover,
\begin{equation*}
y_t=Z-ay,\qquad
y_{tt}=Z_t-aZ+a^2y.
\end{equation*}
Using \eqref{MGT_factorization}, we obtain
\begin{align}\label{MGT_wave_equation}
L_0^{(b)}Z
=\mathcal L_\alpha y-\gamma Z_t+a\gamma Z-a^2\gamma y.
\end{align}
Theorem \ref{CW2} and the uniform bound for $\gamma$ give
\begin{multline}\label{MGT_intermediate}
\mathcal J_s(Z)
\leq C\int_Qe^{2s\varphi_\lambda}
 |\mathcal L_\alpha y|^2\,dxdt
+C\int_Qe^{2s\varphi_\lambda}
 \big(|Z_t|^2+|Z|^2+|y|^2\big)\,dxdt\\
+Cs\int_{\Sigma_0}e^{2s\varphi_\lambda}
 |\partial_nZ|^2\,d\sigma dt.
\end{multline}
Since $P_a y=Z$, $y(\cdot,0)=0$, and $a>0$, Lemma
\ref{lemma_ODE}, together with \eqref{MGT_weight_decay}, yields
\begin{align*}
\int_Q|y|^2e^{2s\varphi_\lambda}\,dxdt
\leq\frac{C}{s}
\int_Q|Z|^2e^{2s\varphi_\lambda}\,dxdt.
\end{align*}
The three lower-order integrals in \eqref{MGT_intermediate} are therefore
absorbed by the second and third terms of $\mathcal J_s(Z)$ for $s$
sufficiently large. This proves \eqref{1.133}.
\end{proof}

\subsection{Lipschitz stability for the relaxation coefficient}

The time differentiation in the Bukhgeim--Klibanov argument produces the
normal trace of $P_a\partial_t(u^1-u^2)$. Consequently, the natural data
space is $H^1$ in time for the single measurement \eqref{MGT_measurement}.

\begin{theorem}\label{mainresult}
Assume \eqref{asumb}--\eqref{asumc}, and let
$\alpha^1,\alpha^2\in\mathcal B_M$. Let $u^k$ denote the solution of
\eqref{eqmain2} corresponding to $\alpha=\alpha^k$, for $k=1,2$, with the same source
and initial data. Suppose that the reference solution satisfies
\begin{align}\label{regsol}
R:=u_{tt}^2\in H^1(0,T;L^\infty(\Omega)),
\qquad
\|R\|_{H^1(0,T;L^\infty(\Omega))}\leq K,
\end{align}
and assume the non-degeneracy condition
\begin{align}\label{mineta}
|u_2(x)|\geq\eta
\qquad\text{for a.e. }x\in\Omega
\end{align}
for some $\eta>0$. Then there exists
\begin{equation*}
C=C(\Omega,T,\Gamma_0,b,c,M,K,\eta)>0
\end{equation*}
such that
\begin{align}\label{stabi}
C^{-1}\|\alpha^1-\alpha^2\|_{L^2(\Omega)}
\leq
\left\|\partial_t\big(\Lambda(\alpha^1)
-\Lambda(\alpha^2)\big)\right\|_{L^2(\Sigma_0)}
\leq
C\|\alpha^1-\alpha^2\|_{L^2(\Omega)}.
\end{align}
In particular,
\begin{align*}\label{stabi_H1}
\|\alpha^1-\alpha^2\|_{L^2(\Omega)}
\leq C\|\Lambda(\alpha^1)-\Lambda(\alpha^2)\|
_{H^1(0,T;L^2(\Gamma_0))}.
\end{align*}
\end{theorem}

\begin{remark}
Condition \eqref{mineta} is an excitation condition on the prescribed
initial acceleration. It ensures that the initial trace generated by a
change in $\alpha$ does not vanish. If additional compatibility conditions
force $u_2$ to vanish at the boundary, a localized version follows exactly
as in Theorem \ref{mainresult2}: one assumes non-degeneracy on
$\Omega_0\Subset\Omega$ and takes the two coefficients equal in
$\Omega\setminus\Omega_0$.
\end{remark}

\begin{remark}
The derivative in \eqref{stabi} does not represent a second experiment.
It is the time derivative of the same boundary trace
$\Lambda(\alpha)$. Its appearance is intrinsic to the proof because the
coefficient is encoded in the initial value of
$(P_a\partial_t(u^1-u^2))_t$.
\end{remark}

\begin{proof}[Proof of Theorem \ref{mainresult}]
Set
\begin{equation*}
z:=u^1-u^2,\qquad
q:=\alpha^2-\alpha^1.
\end{equation*}
Since
$\gamma_{\alpha^2}-\gamma_{\alpha^1}=\alpha^2-\alpha^1=q$,
subtracting the two equations gives
\begin{equation}\label{eqDif}
\left\{
\begin{array}{lll}
\mathcal L_{\alpha^1}z=q(x)R(x,t)&\text{in }Q,\\
z=0&\text{on }\Sigma,\\
z(\cdot,0)=z_t(\cdot,0)=z_{tt}(\cdot,0)=0
&\text{in }\Omega.
\end{array}
\right.
\end{equation}
Let $y:=z_t$. Since $q$ is independent of time, differentiation of
\eqref{eqDif} yields
\begin{equation}\label{eqDer}
\left\{
\begin{array}{lll}
\mathcal L_{\alpha^1}y=q(x)R_t(x,t)&\text{in }Q,\\
y=0&\text{on }\Sigma,\\
y(\cdot,0)=y_t(\cdot,0)=0&\text{in }\Omega.
\end{array}
\right.
\end{equation}
Evaluating \eqref{eqDif} at $t=0$ gives the additional identity
\begin{align}\label{MGT_initial_identity}
(P_a y)_t(x,0)=y_{tt}(x,0)=q(x)u_2(x).
\end{align}

Fix $\lambda\geq\lambda_0$ and choose $s$ sufficiently large. Applying
Theorem \ref{observability} to \eqref{eqDer} and using
\eqref{MGT_initial_identity}, we find
\begin{multline}\label{MGT_weighted_recovery}
\sqrt{s}\int_\Omega e^{2s\varphi_\lambda(x,0)}
 |q(x)u_2(x)|^2\,dx
\leq C\int_Qe^{2s\varphi_\lambda}
 |q(x)R_t(x,t)|^2\,dxdt\\
+Cs\int_{\Sigma_0}e^{2s\varphi_\lambda}
 |\partial_n(P_a y)|^2\,d\sigma dt.
\end{multline}
Since $\varphi_\lambda(x,t)\leq\varphi_\lambda(x,0)$ for
$t\in[0,T]$, assumption \eqref{regsol} implies
\begin{align}\label{MGT_source_bound}
\int_Qe^{2s\varphi_\lambda}|qR_t|^2\,dxdt
\leq K^2\int_\Omega e^{2s\varphi_\lambda(x,0)}
 |q(x)|^2\,dx.
\end{align}
By \eqref{mineta}, the left-hand side of
\eqref{MGT_weighted_recovery} is bounded from below by
\begin{equation*}
\eta^2\sqrt{s}\int_\Omega e^{2s\varphi_\lambda(x,0)}
 |q(x)|^2\,dx.
\end{equation*}
We may therefore absorb \eqref{MGT_source_bound} by choosing $s$
sufficiently large. Once $s$ and $\lambda$ are fixed, the weights are
bounded above and below as needed, and we obtain
\begin{align*}
\|q\|_{L^2(\Omega)}
\leq C\|\partial_n(P_a y)\|_{L^2(\Sigma_0)}.
\end{align*}
Finally, $P_a$ commutes with $\partial_t$, and hence
\begin{align}\label{MGT_trace_identity}
\partial_n(P_a y)
=\partial_n(P_a z_t)
=\partial_t\big(\Lambda(\alpha^1)-\Lambda(\alpha^2)\big).
\end{align}
This proves the first inequality in \eqref{stabi}.

For completeness, we prove the direct estimate as well. Put
$Z:=P_a y$. Equations \eqref{eqDer} and
\eqref{MGT_wave_equation} give the coupled system
\begin{align}\label{MGT_coupled_direct}
\left\{
\begin{aligned}
Z_{tt}-b\Delta Z+\gamma_{\alpha^1}Z_t
-a\gamma_{\alpha^1}Z+a^2\gamma_{\alpha^1}y
&=qR_t,\\
y_t+ay&=Z,\\
Z(\cdot,0)=0,\qquad Z_t(\cdot,0)&=qu_2.
\end{aligned}
\right.
\end{align}
The standard wave-energy estimate, the ODE identity in the second line,
and Gronwall's lemma imply
\begin{align*}\label{MGT_energy_direct}
\|Z\|_{C([0,T];H_0^1(\Omega))}
+\|Z_t\|_{C([0,T];L^2(\Omega))}
+\|y\|_{C([0,T];L^2(\Omega))}
\leq C\|q\|_{L^2(\Omega)}.
\end{align*}
Indeed, $qu_2\in L^2(\Omega)$ and
$qR_t\in L^2(Q)$ by \eqref{regsol}. The hidden-regularity estimate for
the Dirichlet wave equation, applied to the first equation in
\eqref{MGT_coupled_direct}, then yields
\begin{align*}
\|\partial_nZ\|_{L^2(\Sigma)}
\leq C\|q\|_{L^2(\Omega)}.
\end{align*}
Using \eqref{MGT_trace_identity} and restricting the trace to
$\Sigma_0$ proves the second inequality in \eqref{stabi}.
\end{proof}

\begin{corollary}\label{MGT_uniqueness}
Under the assumptions of Theorem \ref{mainresult}, if
\begin{equation*}
\Lambda(\alpha^1)=\Lambda(\alpha^2)
\qquad\text{on }\Sigma_0,
\end{equation*}
then $\alpha^1=\alpha^2$ a.e. in $\Omega$.
\end{corollary}

\section{Concluding remarks and open problems}\label{open}

In this paper, we have established Lipschitz stability results for two
coefficient inverse problems arising in evolution equations with relaxation.
For the Kelvin--Voigt equation, the space-dependent friction coefficient is
recovered from an interior dynamic observation combined with a Cauchy
snapshot at a prescribed time. For the Moore--Gibson--Thompson equation, the
relaxation coefficient is recovered from a single boundary trace under a
standard geometric condition. Although the two equations have different
principal dynamics, the proofs rely on the same structural principle: an
exact PDE--ODE factorization that allows classical parabolic and hyperbolic
Carleman estimates to be combined with weighted estimates for a first-order
equation in time.

The analysis also raises several questions that remain open. We briefly
discuss three of them.

\begin{itemize}
\item \emph{Space--time-dependent coefficients.}
It is natural to ask whether the stability results remain valid when
$b=b(x,t)$ or $\alpha=\alpha(x,t)$. This is not a direct extension of the
arguments developed here. Indeed, differentiating the difference equation
in time introduces the additional unknowns $\partial_t(b^1-b^2)$ and
$\partial_t(\alpha^1-\alpha^2)$, while the fixed-time identities recover
only one temporal slice of the coefficient. A plausible intermediate problem
is to consider a prescribed temporal profile, for example
\begin{equation*}
b(x,t)=b_\ast(x,t)+\rho(t)\beta(x),
\qquad
\alpha(x,t)=\alpha_\ast(x,t)+\mu(t)\kappa(x),
\end{equation*}
where $\rho$ and $\mu$ are known and the spatial amplitudes $\beta$ and
$\kappa$ are unknown. Recovering a general space--time coefficient is
expected to require additional observations or a different Carleman
strategy.

\item \emph{Variable principal coefficients in the MGT equation.}
Another question is whether the diffusivity $b$ or the sound speed $c$ can
depend on space or time and, possibly, be recovered from measurements. Such
dependence changes the principal part of the equation and produces
commutator terms in the wave--ODE decomposition, since $c^2/b$ is no longer
constant. It would be interesting to determine which factorizations and
geometric conditions survive in this setting, as well as the minimal number
of experiments needed to identify more than one coefficient.

\item \emph{Numerical reconstruction.}
The Lipschitz estimates proved in this paper provide a natural starting point
for reconstructing $b(x)$ and $\alpha(x)$. A further problem is to design and
analyze a stable numerical procedure, for instance through a
Carleman-weighted functional or an iterative regularization method, and to
study its behavior under discretization and noisy observations.
\end{itemize}


\end{document}